# THE GEOMETRY OF FOURTH ORDER DIFFERENTIAL OPERATOR

ROHOLLAH BAKHSHANDEH-CHAMAZKOTI

ABSTRACT. In present paper, the equivalence problem for fourth order differential operators with one variable under general fiber-preserving transformation using the Cartan method of equivalence is applied. Two versions of equivalence problems are considered. First, the direct equivalence problem and second equivalence problem is to determine the sufficient and necessary conditions on two fourth order differential operators such that there exists a fiber-preserving transformation mapping one to the other according to gauge equivalence.

## Contents

## 1. Introduction

The classification of linear differential equations is a special case of the general problem of classifying differential operators, which has a variety of important applications, including quantum mechanics and the projective geometry of curves [1]. In this attempt we shall solve the method of local equivalence problem by three versions of the equivalence problem for the class of linear fourth order operators on the line. For simplicity, we shall only deal with the local equivalence problem for scalar differential operators in a single independent variable, although these problems are important for matrix-valued and partial differential operators as well.

The general equivalence problem is to recognize when two geometrical objects are mapped on each other by a certain class of diffeomorphisms. E. Cartan developed the general equivalence problem and provided a systematic procedure for determining the necessary and sufficient condition [2, 3]. In Cartan's approach, the conditions of equivalence of two objects must be reformulated in terms of differential forms. We associate a collection of one-forms to an object under investigation in the original coordinates; the corresponding object in the new coordinates will have its own collection of one-forms. Once an equivalence problem has been reformulated in the proper Cartan form, in terms of a coframe $\omega$ on the $m$-dimensional base manifold $M$, along with a structure group $G \subset \text{GL}(m)$, we can apply the Cartan equivalence method. The goal is to normalize the structure group valued coefficients in a suitably invariant manner, and this is accomplished through the determination of a sufficient number of invariant combinations thereof [1].

The problems here are related to the more general equivalence problem for fourth order ordinary differential equations which E. Cartan studied under point transformations [4], and S. S. Chern turned his attention to the problem under contact transformations [5] and Hajime Sato et all [6], but are specialized by linearity.

Niky Kamran and Peter J. Olver have been solved equivalence problem for second order differential operator with two versions of the equivalence problem [7] and also Nadjafikhah and Bakhshandeh have been solved this problem for fourth order operators [8]. They didn't do the projective case because all (nonsingular) second order differential operators are projectively equivalent, and so the second order case is not interesting. But we also solve the full projective equivalence problem for fourth order differential operators. Projective problems was discussed at length in [9] and it also has implications for equivalence problems for curves in projective space. A brief survey of Wilczynski's analysis can be found starting on [1]. Extensions of Wilczynski's work to nonlinear ordinary differential equations can be found in the paper [10].

## 2. Equivalence of fourth order differential operators

Consider the fourth order differential operator applied on a scalar-valued function $u(x)$

$$(2.1) \qquad \mathcal{D}[u] = \sum_{i=0}^{4} f_i(x) \, D^i u$$

and another fourth order differential operator applied on a scalar-valued function $\bar{u}(\bar{x})$

$$(2.2) \qquad \bar{\mathcal{D}}[\bar{u}] = \sum_{i=0}^{4} \bar{f}_i(\bar{x}) \, \bar{D}^i \bar{u}.$$

where $f_i$ and $\bar{f}_i$, $i = 1, 2, 3, 4$, are analytic functions of the real variable $x$ and $\bar{x}$ respectively. Further, $D^i = d/dx^i$, $\bar{D}^i = d/d\bar{x}^i$ and $D^0 = \bar{D}^0 = \mathrm{Id}$ are the identity operators.

The appropriate space to work in will be the fourth jet space $\mathrm{J}^4$, which has local coordinates

$$\Upsilon = \{(x, u, p, q, r, s) \in \mathrm{J}^4 : p = u_x, q = u_{xx}, r = u_{xxx}, s = u_{xxxx}\}$$

and our goal is to know whether there exists a suitable transformation of variables $(x, u, p, q, r, s) \longrightarrow (\bar{x}, \bar{u}, \bar{p}, \bar{q}, \bar{r}, \bar{s})$ which brings (2.1) to (2.2). Several types of such transformations are of particular importance. Here we consider fiber preserving transformations, which are of the form

$$(2.3) \qquad \bar{x} = \xi(x), \qquad \bar{u} = \varphi(x) \, u,$$

where $\varphi(x) \neq 0$. Using the chain rule formula we find following relation between the total derivative operators

$$(2.4) \qquad \bar{D} = \frac{d}{d\bar{x}} = \frac{1}{\xi'(x)} \frac{d}{dx} = \frac{1}{\xi'(x)} \, D.$$

We first consider the *direct equivalence problem*, which identifies the two linear differential functions

$$(2.5) \qquad \mathcal{D}[u] = \bar{\mathcal{D}}[\bar{u}].$$

under change of variables (2.3). This induces the transformation rule

$$\bar{\mathcal{D}} = \mathcal{D} \cdot \frac{1}{\varphi(x)} \qquad \text{when} \qquad \bar{x} = \xi(x), \tag{2.6}$$

on the differential operators themselves, and solving local direct equivalence problem is to find explicit conditions on the coefficients of the two differential operators that guarantee that they satisfy (2.5) for some change of variables of the form (2.3).

The transformation rule (2.6) doesn't preserve either the eigenvalue problem $\mathcal{D}[u] = \lambda u$ or the Schrödinger equation $iu_t = \mathcal{D}[u]$, since we are missing a factor of $\varphi(x)$. For solving this problem, we consider the *gauge equivalence* with the following transformation rule

$$\bar{\mathcal{D}} = \varphi(x) \cdot \mathcal{D} \cdot \frac{1}{\varphi(x)} \qquad \text{when} \qquad \bar{x} = \xi(x). \tag{2.7}$$

**Proposition 1.** *Suppose $\mathcal{D}$ and $\bar{\mathcal{D}}$ be fourth-order differential operators. There are two coframes $\Omega = \{\omega^1, \omega^2, \omega^3, \omega^4, \omega^5, \omega^6\}$ and $\bar{\Omega} = \{\bar{\omega}^1, \bar{\omega}^2, \bar{\omega}^3, \bar{\omega}^4, \bar{\omega}^5, \bar{\omega}^6\}$ on open subsets $\Gamma$ and $\bar{\Gamma}$ of the fourth jet space, respectively, such that the differential operators are equivalent under the pseudogroup (2.3) according to the respective transformation rules (2.6) and (2.7) which coframes $\Omega$ and $\bar{\Omega}$ satisfy in following relation*

$$\begin{pmatrix} \bar{\omega}^1 \\ \bar{\omega}^2 \\ \bar{\omega}^3 \\ \bar{\omega}^4 \\ \bar{\omega}^5 \\ \bar{\omega}^6 \end{pmatrix} = \begin{pmatrix} a_1 & 0 & 0 & 0 & 0 & 0 \\ 0 & 1 & 0 & 0 & 0 & 0 \\ 0 & a_2 & a_3 & 0 & 0 & 0 \\ 0 & a_4 & a_5 & a_6 & 0 & 0 \\ 0 & a_7 & a_8 & a_9 & a_{10} & 0 \\ 0 & 0 & 0 & 0 & 0 & 1 \end{pmatrix} \begin{pmatrix} \omega^1 \\ \omega^2 \\ \omega^3 \\ \omega^4 \\ \omega^5 \\ \omega^6 \end{pmatrix} \tag{2.8}$$

*where $a_i \in \mathbb{R}$ for $i = 1, \cdots, 10$ and $a_1 a_3 a_6 a_{10} \neq 0$.*

*Proof.* Note first that a point transformation will be in the desired linear form (2.3) if and only if, for pair of functions $\alpha = \xi_x$ and $\beta = \varphi_x/\varphi$, one-form equations

$$d\bar{x} = \alpha\, dx, \tag{2.9}$$

$$\frac{d\bar{u}}{\bar{u}} = \frac{du}{u} + \beta\, dx. \tag{2.10}$$

hold on the subset of $J^4$ where $u \neq 0$. In order that the derivative variables $p, q, r$ and $s$ transform correctly, we need to preserve the contact ideal $\mathcal{I}$ on $J^4$, which is

$$\mathcal{I} = \langle du - p\, dx, dp - q\, dx, dq - r\, dx, dr - s\, dx \rangle. \tag{2.11}$$

Generally, a diffeomorphism $\Phi : J^4 \to J^4$ determines a contact transformation if and only if

$$d\bar{u} - \bar{p}\, d\bar{x} = a_1(du - p\, dx), \tag{2.12}$$

$$d\bar{p} - \bar{q}\, d\bar{x} = a_2(du - p\, dx) + a_3(dp - q\, dx), \tag{2.13}$$

$$d\bar{q} - \bar{r}\, d\bar{x} = a_4(du - p\, dx) + a_5(dp - q\, dx) + a_6(dq - r\, dx), \tag{2.14}$$

$$d\bar{r} - \bar{s}\, d\bar{x} = a_7(du - p\, dx) + a_8(dp - q\, dx) + a_9(dq - r\, dx) + a_{10}(dr - s\, dx), \tag{2.15}$$

where $a_i$ are functions on $J^4$. The combination of the first contact condition (2.12) with the linearity conditions (2.9) and (2.10) constitutes part of an overdetermined equivalence problem.

Taking $\beta = -p/u$, $a_1 = 1/u$, in (2.10) and (2.12), it is found the one-form

$$\frac{d\bar{u} - \bar{p}\, d\bar{x}}{\bar{u}} = \frac{du - p\, dx}{u}, \tag{2.16}$$

which is invariant, and (2.16) can replace both (2.10) and (2.12). Therefore, we may choose five elements of our coframe the one-forms

$$\omega^1 = dx, \quad \omega^2 = \frac{du - p\, dx}{u}, \quad \omega^3 = dp - q\, dx, \quad \omega^4 = dq - r\, dx, \quad \omega^5 = dr - s\, dx, \tag{2.17}$$

which are defined on the fourth jet space $J^4$ locally parameterized by $(x, u, p, q, r, s)$, with the transformation rules

$$\bar{\omega}^1 = a_1 \omega^1, \quad \bar{\omega}^2 = \omega^2, \quad \bar{\omega}^3 = a_2 \omega^2 + a_3 \omega^3, \quad \bar{\omega}^4 = a_4 \omega^2 + a_5 \omega^3 + a_6 \omega^4 \tag{2.18}$$
$$\bar{\omega}^5 = a_7 \omega^2 + a_8 \omega^3 + a_9 \omega^4 + a_{10} \omega^5.$$

According to (2.5), the function $I(x, u, p, q, r, s) = \mathcal{D}[u] = f_4(x)s + f_3(x)r + f_2(x)q + f_1(x)p + f_0(x)u$ is an invariant for the problem, and thus its differential

$$\omega^6 = dI = f_4 ds + f_3 dr + f_2 dq + f_1 dp + f_0 du + (f_4' s + f_3' r + f_2' q + f_1' p + f_0' u) dx, \tag{2.19}$$

is an invariant one-form, thus one can take it as a final element of our coframe.

In the second problem (2.7), for the extra factor of $\varphi$, the invariant is

$$I(x, u, p, q, r, s) = \frac{\mathcal{D}[u]}{u} = \frac{f_4(x)ds + f_3(x)r + f_2(x)q + f_1(x)p}{u} + f_0(x). \tag{2.20}$$

Thus, it is found

$$\omega^6 = dI = \frac{f_4}{u} ds + \frac{f_3}{u} dr + \frac{f_2}{u} dq + \frac{f_1}{u} dp - \frac{f_3 r + f_2 q + f_1 p}{u^2} du + \left\{ \frac{f_3' r + f_2' q + f_1' p}{u} + f_0' \right\} dx, \tag{2.21}$$

as a final element of coframe for the equivalence problem (2.7). The set of one-forms

$$\Omega = \{\omega^1, \omega^2, \omega^3, \omega^4, \omega^5, \omega^6\}$$

is a coframe on the subset

$$\Gamma^* = \left\{ (x, u, p, q, r, s) \in J^4 \,\Big|\, u \neq 0 \text{ and } f_4(x) \neq 0 \right\}. \tag{2.22}$$

All of attention is restricted to a connected component $\Gamma \subset \Gamma^*$ of the subset (2.22) that the signs of $f_0(x)$ and $u$ are fixed. It means that the last coframe elements agree up to contact

$$\bar{\omega}^6 = \omega^6. \tag{2.23}$$

Viewing (2.18) and (2.23) relations, one can find the structure group associated with the equivalence problems (2.6) and (2.7) that is a ten-dimensional matrix group $G$ such that $\bar{\Omega} = G\Omega$ which leads to (2.8) and then the *lifted coframe* on the space $J^4 \times G$ has the form

$$\begin{aligned}
\theta^1 &= a_1 \omega^1, \\
\theta^2 &= \omega^2, \\
\theta^3 &= a_2 \omega^2 + a_3 \omega^3, \\
\theta^4 &= a_4 \omega^2 + a_5 \omega^3 + a_6 \omega^4, \\
\theta^5 &= a_7 \omega^2 + a_8 \omega^3 + a_9 \omega^4 + a_{10} \omega^5, \\
\theta^5 &= \omega^6.
\end{aligned} \tag{2.24}$$

☐

Here, the main results are presented as following two theorems:

**Theorem 1.** *The final structure equations for direct equivalence with (2.17) and (2.19) coframes are*

(2.25)
$$d\theta^1 = \frac{1}{4}\theta^1 \wedge \theta^2,$$
$$d\theta^2 = \theta^1 \wedge \theta^3,$$
$$d\theta^3 = \theta^1 \wedge \theta^4 + \frac{1}{4}\theta^2 \wedge \theta^3,$$
$$d\theta^4 = I_1\theta^1 \wedge \theta^3 + \theta^1 \wedge \theta^5 + \frac{1}{2}\theta^2 \wedge \theta^4,$$
$$d\theta^5 = I\theta^1 \wedge \theta^2 + I_2\theta^1 \wedge \theta^3 + I_3\theta^1 \wedge \theta^4 + \theta^1 \wedge \theta^6 + \frac{3}{4}\theta^2 \wedge \theta^5 + 3\theta^3 \wedge \theta^4,$$
$$d\theta^6 = 0,$$

where the coefficients $I_1, I_2, I_3$ and $I$ are

(2.26)
$$I_1 = -\frac{\sqrt[4]{f_4 u}}{2f_4 u}\left[5f_4 p - 3\dot{f}_4 u + 2f_3 u\right]$$
$$I_2 = \frac{1}{64 f_4^2 u^2 \sqrt[4]{f_4 u}}\Big[(20 f_3 \dot{f}_4^2 - 16 f_3 f_4 \ddot{f}_4 - 45 \dot{f}_4^3 - 16 f_4^2 \dddot{f}_4 + 64 f_1 f_4^2 - 16 f_2 f_4 \dot{f}_4 + 60 f_4 \dot{f}_4 \ddot{f}_4)u^3$$
$$+ 35 f_4^3 p^3 + 240 f_4 u^2 r + (65 f_4 \dot{f}_4^2 - 52 f_4^2 \ddot{f}_4 + 112 f_2 f_4^2 - 40 f_3 f_4 \dot{f}_4)u^2 p - 100 f_4^3 upq$$
$$+ (17\dot{f}_4 - 28 f_3) f_4^2 up^2 + (176 f_3 \dot{f}_4^2 - 84 f_4^2 \dot{f}_4)u^2 q\Big]$$
$$I_3 = -\frac{1}{16 f_4 \sqrt{f_4 uu}}\Big[(16 f_2 f_4 + 5\dot{f}_4^2 - 16 f_4 \dot{f}_3 + 8 f_4 \ddot{f}_4)u^2 + 5 f_4 p^2 + 40 f_4^2 uq$$
$$+ (32 f_3 - 38 \dot{f}_4) f_4 up\Big],$$
$$I = -(f_4 s + f_3 r + f_2 q + f_1 p + f_0 u).$$

**Theorem 2.** *The final structure equations for gauge equivalence with (2.17) and (2.21) coframes are*

(2.27)
$$d\theta^1 = 0,$$
$$d\theta^2 = \theta^1 \wedge \theta^3,$$
$$d\theta^3 = \theta^1 \wedge \theta^4,$$
$$d\theta^4 = I_1\theta^1 \wedge \theta^3 + I_2\theta^1 \wedge \theta^4 + \theta^1 \wedge \theta^5,$$
$$d\theta^5 = I_3\theta^1 \wedge \theta^3 + I_4\theta^1 \wedge \theta^4 + \theta^1 \wedge \theta^6,$$
$$d\theta^6 = 0,$$

where the coefficients $I_1, \ldots, I_4$ are

$$
\begin{aligned}
I_1 &= -\frac{1}{4f_4\sqrt{f_4}u^2}\left[16f_4^2 uq + 8f_4^2 p^2 + (8f_3f_4 - 10f_4\dot{f}_4)up + (2\dot{f}_4^2 - f_4\ddot{f}_4 - f_3\dot{f}_4)u^2\right], \\
I_2 &= -\frac{1}{2\sqrt[4]{f_4^3}u}\left[8f_4 p + (2f_3 - 3\dot{f}_4)u\right], \\
I_3 &= \frac{1}{f_4^2\sqrt[4]{f_4}u^3}\Big[(64f_3f_4^2 - 64f_4^2\dot{f}_4)u^2 q + 256f_4^3 p^3 + (128f_3f_4^2 - 244f_4^2\dot{f}_4)up^2 \\
&\quad + (128f_2f_4^2 - 128f_4^2\dot{f}_3 + 96f_4^2\ddot{f}_4 + 16f_3f_4\dot{f}_4)u^2 p + (16f_4\dot{f}_3\dot{f}_4 - 16f_2 f_4\dot{f}_4 \\
&\quad - 12f_4\dot{f}_4\ddot{f}_4 - 4f_3\dot{f}_4^2 + 3\dot{f}_4^3 + 64f_1 f_4^2)u^3\Big], \\
I_4 &= -\frac{1}{16f_4\sqrt{f_4}u^2}\Big[32f_4^2 p^2 - 32f_4^2 uq + (16f_3f_4 - 24f_4\dot{f}_4)up + (12f_4\ddot{f}_4 - 16f_4\dot{f}_3 + 4f_3\dot{f}_4 \\
&\quad - 3\dot{f}_4^2 + 16f_2 f_4)u^2\Big]
\end{aligned}
$$
(2.28)

## 3. The proof of Theorem 1

First, the initial five one-forms (2.17) and (2.19) are taken as our final coframe constituent. The next step is to calculate the differentials of lifted coframe elements (2.24). An explicit computation leads to the structure equations

$$
\begin{aligned}
d\theta^1 &= \alpha^1 \wedge \theta^1, \\
d\theta^2 &= T_{12}^2 \theta^1 \wedge \theta^2 + T_{13}^2 \theta^1 \wedge \theta^3, \\
d\theta^3 &= \alpha^2 \wedge \theta^2 + \alpha^3 \wedge \theta^3 + T_{12}^3 \theta^1 \wedge \theta^2 + T_{13}^3 \theta^1 \wedge \theta^3 + T_{14}^3 \theta^1 \wedge \theta^4, \\
d\theta^4 &= \alpha^4 \wedge \theta^2 + \alpha^5 \wedge \theta^3 + \alpha^6 \wedge \theta^4 + T_{12}^4 \theta^1 \wedge \theta^2 + T_{13}^4 \theta^1 \wedge \theta^3 \\
&\quad + T_{14}^4 \theta^1 \wedge \theta^4 + T_{15}^4 \theta^1 \wedge \theta^5, \\
d\theta^5 &= \alpha^7 \wedge \theta^2 + \alpha^8 \wedge \theta^3 + \alpha^9 \wedge \theta^4 + \alpha^{10} \wedge \theta^5 + T_{12}^4 \theta^1 \wedge \theta^2 + T_{13}^4 \theta^1 \wedge \theta^3 \\
&\quad + T_{14}^4 \theta^1 \wedge \theta^4 + T_{15}^4 \theta^1 \wedge \theta^5 + T_{16}^5 \theta^1 \wedge \theta^6, \\
d\theta^6 &= 0,
\end{aligned}
$$
(3.1)

with

$$
\alpha^1 = \frac{da_1}{a_1}, \quad \alpha^2 = \frac{a_3 da_2 - a_2 da_3}{a_3}, \quad \alpha^3 = \frac{da_3}{a_3},
$$

$$
\alpha^4 = \frac{a_3 a_6 da_4 - a_2 a_6 da_5 + (a_2 a_5 - a_3 a_4)da_6}{a_3 a_6}, \quad \alpha^5 = \frac{a_6 da_5 - a_5 da_6}{a_3 a_6}, \quad \alpha^6 = \frac{da_6}{a_6},
$$

$$
\alpha^7 = \frac{a_3 a_6 a_{10} da_7 - a_6 a_{10} da_8 + a_{10}(a_2 a_5 - a_3 a_4)da_9 - (a_3 a_6 a_7 - a_3 a_4 a_9 - a_2 a_6 a_8 + a_2 a_5 a_9)da_{10}}{a_3 a_6 a_{10}},
$$

$$
\alpha^8 = \frac{a_6 a_{10} da_8 - a_5 a_{10} da_9 + (a_5 a_9 - a_6 a_8)da_{10}}{a_3 a_6 a_{10}}, \quad \alpha^9 = \frac{a_{10} da_9 - a_9 da_{10}}{a_6 a_{10}}, \quad \alpha^{10} = \frac{da_{10}}{a_{10}},
$$

forming a basis for the right-invariant *Maurer-Cartan forms* on the Lie group $G$. The essential torsion coefficients are

(3.2) $\quad T_{12}^2 = -\frac{a_2 + a_3 p}{a_1 a_3 u}, \quad T_{13}^2 = \frac{1}{a_1 a_3 u}, \quad T_{14}^3 = \frac{a_3}{a_1 a_6}, \quad T_{15}^4 = \frac{a_6}{a_1 a_{10}}, \quad T_{16}^5 = \frac{a_{10}}{a_1 f_4}.$

It is possible to normalize the group parameters by setting

(3.3) $\quad a_1 = \frac{1}{\sqrt[4]{f_4 u}}, \quad a_2 = -\frac{\sqrt[4]{f_4 u}}{u}p, \quad a_3 = \frac{\sqrt[4]{f_4 u}}{u}, \quad a_6 = \frac{\sqrt{f_4 u}}{u}, \quad a_{10} = -\frac{f_4}{\sqrt[4]{f_4 u}}.$

In the second loop, the normalization (3.3) is substituted in the lifted coframe (2.24) and calculate the differentials of new invariant coframe to obtain revised structure equations. Now, the essential torsion components (3.2) are normalized by the parameters

$$(3.4) \qquad a_5 = \frac{\dot{f}_4 u - 7 f_4 p}{4 \sqrt[4]{(f_4 u)^3}} \qquad a_9 = \frac{(4 f_3 - 3 \dot{f}_4) u + f_4 p}{4 \sqrt[4]{(f_4 u)^3}}.$$

To determine the remaining parameters $a_4, a_7, a_8$, the obtained parameters (3.4) are substituted into (2.24), and recalculate the differentials. Therefore, the new structure equations are

$$(3.5) \quad \begin{aligned} d\theta^1 &= \frac{1}{4} \theta^1 \wedge \theta^2, \\ d\theta^2 &= \theta^1 \wedge \theta^3, \\ d\theta^3 &= T^3_{12} \theta^1 \wedge \theta^2 + \frac{1}{4} \theta^2 \wedge \theta^3 + \theta^1 \wedge \theta^4, \\ d\theta^4 &= \alpha^4 \wedge \theta^2 + T^4_{12} \theta^1 \wedge \theta^2 + T^4_{13} \theta^1 \wedge \theta^3 + T^4_{14} \theta^1 \wedge \theta^4 + T^4_{23} \theta^2 \wedge \theta^3 + \theta^1 \wedge \theta^5 + \frac{1}{2} \theta^2 \wedge \theta^4, \\ d\theta^5 &= \alpha^8 \wedge \theta^3 + \alpha^7 \wedge \theta^2 + T^5_{12} \theta^1 \wedge \theta^2 + T^5_{13} \theta^1 \wedge \theta^3 + T^5_{23} \theta^2 \wedge \theta^3 + T^5_{14} \theta^1 \wedge \theta^4 \\ &\quad - \frac{1}{4} \theta^2 \wedge \theta^5 + \frac{1}{4} \theta^3 \wedge \theta^4 - \frac{1}{4} \theta^2 \wedge \theta^5 + \theta^1 \wedge \theta^6, \\ d\theta^6 &= 0. \end{aligned}$$

where $\alpha^4, \alpha^7$ and $\alpha^8$ are the Maurer-Cartan forms on $G$ and the essential torsion coefficients are

$$(3.6) \qquad T^3_{12} = -\frac{a_4 f_4 u \sqrt{f_4 u} + f_4^2 u q}{f_4 u \sqrt{f_4 u}}.$$

By assumption $f_4 u \neq 0$ then one can do following normalization by setting

$$(3.7) \qquad a_4 = -\frac{f_4 q}{\sqrt{f_4 u}}.$$

Substituting (3.7) in (2.24) and recomputing the differentials leads to

$$(3.8) \quad \begin{aligned} d\theta^1 &= \frac{1}{4} \theta^1 \wedge \theta^2, \\ d\theta^2 &= \theta^1 \wedge \theta^3, \\ d\theta^3 &= \frac{1}{4} \theta^2 \wedge \theta^3 + \theta^1 \wedge \theta^4, \\ d\theta^4 &= T^4_{12} \theta^1 \wedge \theta^2 + T^4_{13} \theta^1 \wedge \theta^3 + T^4_{14} \theta^1 \wedge \theta^4 + \theta^1 \wedge \theta^5 + \frac{1}{2} \theta^2 \wedge \theta^4, \\ d\theta^5 &= \alpha^8 \wedge \theta^3 + \alpha^7 \wedge \theta^2 + T^5_{12} \theta^1 \wedge \theta^2 + T^5_{13} \theta^1 \wedge \theta^3 + T^5_{23} \theta^2 \wedge \theta^3 + T^5_{14} \theta^1 \wedge \theta^4 \\ &\quad - \frac{1}{4} \theta^2 \wedge \theta^5 + \frac{1}{4} \theta^3 \wedge \theta^4 - \frac{1}{4} \theta^2 \wedge \theta^5 + \theta^1 \wedge \theta^6, \\ d\theta^6 &= 0. \end{aligned}$$

Finally, with a simple calculation one can find remained normalization by

$$(3.9) \quad \begin{aligned} a_7 &= -\frac{1}{\sqrt[4]{f_4 u u}} \left[ f_4 p q + 4 f_4 u q + (4 f_3 - 3 \dot{f}_4) u q \right], \\ a_8 &= -\frac{\sqrt{f_4 u}}{4 f_4^2 u} \left[ 11 f_4^2 q + (2 \dot{f}_4^2 - f_4 \ddot{f}_4 - f_3 \dot{f}_4) u + (7 f_3 f_4 - 8 f_4 \dot{f}_4) p \right], \end{aligned}$$

and then it leads to final structure equations (2.27) with coefficients (2.26).

## 4. The proof of Theorem 2

The Cartan formulation of gauge equivalence of fourth differential operators will use the same initial five coframes (2.17), but now the final 1-form element is (2.21). In the first loop through the second equivalence problem procedure, according to Proposition 1, the structure group $G$ in (2.8) relation is exactly the structure group of direct equivalence, and then the equivalence method has the same intrinsic structure (3.1) by the essential torsion coefficients

$$(4.1) \quad T_{12}^2 = -\frac{a_2 + a_3 p}{a_1 a_3 u}, \quad T_{13}^2 = \frac{1}{a_1 a_3 u}, \quad T_{14}^3 = \frac{a_3}{a_1 a_6}, \quad T_{15}^4 = \frac{a_6}{a_1 a_{10}}, \quad T_{15}^5 = \frac{a_{10} u}{a_1 f_4}.$$

One can normalize the group parameters by setting

$$(4.2) \quad a_1 = \frac{1}{\sqrt[4]{f_4}}, \quad a_2 = -\frac{\sqrt[4]{f_4}}{u} p, \quad a_3 = \frac{\sqrt[4]{f_4}}{u}, \quad a_6 = \frac{\sqrt{f_4}}{u}, \quad a_{10} = \frac{\sqrt[4]{f_4^3}}{u}.$$

In the second loop of the present equivalence problem, the normalization (4.2) is substituted in lifted coframe (2.24) and calculate differentials of new invariant coframe for finding following revised structure equations:

$$(4.3) \quad \begin{aligned} d\theta^1 &= 0, \\ d\theta^2 &= \theta^1 \wedge \theta^3, \\ d\theta^3 &= T_{12}^3 \theta^1 \wedge \theta^2 + T_{13}^3 \theta^1 \wedge \theta^3 + \theta^1 \wedge \theta^4, \\ d\theta^4 &= \alpha^4 \wedge \theta^2 + \alpha^5 \wedge \theta^3 + T_{12}^4 \theta^1 \wedge \theta^2 + T_{13}^4 \theta^1 \wedge \theta^3 + T_{14}^4 \theta^1 \wedge \theta^4 + T_{23}^4 \theta^2 \wedge \theta^3 \\ &\quad + \theta^1 \wedge \theta^5 - \theta^2 \wedge \theta^4, \\ d\theta^5 &= \alpha^7 \wedge \theta^2 + \alpha^8 \wedge \theta^3 + \alpha^9 \wedge \theta^4 + T_{12}^5 \theta^1 \wedge \theta^2 + T_{13}^5 \theta^1 \wedge \theta^3 + T_{14}^5 \theta^1 \wedge \theta^4 \\ &\quad + T_{15}^5 \theta^1 \wedge \theta^5 + T_{23}^5 \theta^2 \wedge \theta^3 - \theta^2 \wedge \theta^5 + \theta^1 \wedge \theta^6, \\ d\theta^6 &= 0, \end{aligned}$$

where $\alpha^4, \alpha^5, \alpha^7, \alpha^8$ and $\alpha^9$ are the Maurer-Cartan forms and the essential torsion components of structure equations (4.3) are

$$(4.4) \quad \begin{aligned} T_{12}^3 &= -\frac{a_4 \sqrt{f_4} u + f_4 q}{\sqrt{f_4} u}, \qquad T_{13}^3 = -\frac{8 f_4 p + (4 a_5 \sqrt[4]{f_4^3} u - \dot{f}_4) u}{4 \sqrt[4]{f_4^3} u}, \\ T_{15}^5 &= \frac{(4 a_9 \sqrt[4]{f_4^3} - 4 f_3 + 3 \dot{f}_4) u - 4 f_4 p}{4 \sqrt[4]{f_4^3} u}, \end{aligned}$$

and so the normalization is

$$(4.5) \quad a_4 = -\frac{\sqrt{f_4}}{u} q, \quad a_5 = -\frac{8 f_4 p - \dot{f}_4 u}{4 \sqrt[4]{f_4^3} u}, \quad a_9 = \frac{4 f_4 p + (4 f_3 - 3 \dot{f}_4) u}{4 \sqrt[4]{f_4^3} u}.$$

Putting (4.5) into (2.24) and then recomputing the differential of new 1-forms leads to

$$\begin{aligned}
d\theta^1 &= 0, \\
d\theta^2 &= \theta^1 \wedge \theta^3, \\
d\theta^3 &= \theta^1 \wedge \theta^4, \\
d\theta^4 &= T^4_{12}\theta^1 \wedge \theta^2 + T^4_{13}\theta^1 \wedge \theta^3 + T^4_{14}\theta^1 \wedge \theta^4 + \theta^1 \wedge \theta^5, \\
d\theta^5 &= \alpha^7 \wedge \theta^2 + \alpha^8 \wedge \theta^3 + T^5_{12}\theta^1 \wedge \theta^2 + T^5_{13}\theta^1 \wedge \theta^3 + T^5_{14}\theta^1 \wedge \theta^4 \\
&\quad + T^5_{15}\theta^1 \wedge \theta^5 + T^5_{23}\theta^2 \wedge \theta^3 - \theta^2 \wedge \theta^5 + \theta^3 \wedge \theta^4 + \theta^1 \wedge \theta^6, \\
d\theta^6 &= 0.
\end{aligned} \quad (4.6)$$

This immediately implies following normalization

$$(4.7) \quad a_7 = -\frac{4f_4 pq + 4f_4 ru + (4f_3 - 3\dot{f}_4)qu}{4\sqrt[4]{f_4}\, u^2}, \qquad a_8 = \frac{\sqrt{f_4}}{u}\, q.$$

Thus the final invariant coframe is now given by

$$\begin{aligned}
\theta^1 &= \frac{dx}{\sqrt[4]{f_4}}, \\
\theta^2 &= \frac{du - p\, dx}{u}, \\
\theta^3 &= \frac{\sqrt[4]{f_4}}{u^2}\Big[(p^2 - qu)\, dx - p\, du + u\, dp\Big], \\
\theta^4 &= -\frac{1}{4\sqrt{f_4}\, u^3}\Big[(4f_4 u^2 r + \dot{f}_4 u^2 q - \dot{f}_4 u p^2 - 12 f_4 upq + 8 f_4 p^3)\, dx \\
&\qquad\qquad + (\dot{f}_4 up + 4 f_4 uq - 8 f_4 p^2)\, du + (8 f_4 p - \dot{f}_3 u)u\, dp - 4 f_4 u^2\, dq\Big], \\
\theta^5 &= \frac{1}{4\sqrt{f_4}\, u^3}\Big[(8 f_4 p^2 q - 4 f_4 uq^2 - 4 f_4 u^2 s + 4 f_3 u^2 r (4f_3 - 3\dot{f}_4)upq + 3\dot{f}_4 u^2 r)\, dx \\
&\qquad\qquad + (8 f_4 pq + 4 f_4 ur + (4f_3 - 3\dot{f}_4)uq)\, du + (4 f_4 uq)\, dp + (4 f_4 p + (4f_3 - \dot{f}_4)u)u\, dq \\
&\qquad\qquad - 4 f_4 u^2\, dr\Big], \\
\theta^6 &= \frac{f'_4 s + f'_3 r + f'_2 q + f_1 p + f'_0 u}{u}\, dx - \frac{f_4 r + f_3 r + f_2 q + f_1 p}{u^2}\, du + \frac{f_1}{u}\, dp \\
&\qquad + \frac{f_2}{u}\, dq + \frac{f_3}{u} dr + \frac{f_4}{u} ds.
\end{aligned} \quad (4.8)$$

Then the final structure equations (2.27) with fundamental invariant coefficients (2.28) are obtained.

Department of Mathematics, Faculty of Basic science, Babol University of Technology, Babol, Iran.
*E-mail address*: r_bakhshandeh@nit.ac.ir